\documentclass[12pt]{article}

\usepackage{amsmath, epsfig, cite}
\usepackage{amssymb}
\usepackage{amsfonts}
\usepackage{latexsym}

\newtheorem{thm}{Theorem}[section]

\newtheorem{conj}[thm]{Conjecture}

\numberwithin{equation}{section}

\newcommand{\qed}{{\hfill$\square$}\medskip}

\begin{document}


\begin{center}
{\Large\bf Proof of some conjectures of Z.-W. Sun on\\[5pt] the divisibility of certain double-sums}
\end{center}

\vskip 2mm \centerline{Victor J. W. Guo$^1$ and Ji-Cai Liu$^2$}
\begin{center}
{\footnotesize Department of Mathematics, Shanghai Key Laboratory of
PMMP, East China Normal University,\\ 500 Dongchuan Road, Shanghai
200241,
 People's Republic of China\\
$^1${\tt jwguo@math.ecnu.edu.cn,\quad
http://math.ecnu.edu.cn/\textasciitilde{jwguo}}  \quad {$^2$\tt jc2051@163.com} }
\end{center}


\vskip 0.7cm \noindent{\bf Abstract.} Z.-W. Sun introduced three kinds of numbers:
\begin{align*}S_n=\sum_{k=0}^{n}{n\choose k}^2{2k\choose k}(2k+1),\qquad
s_n=\sum_{k=0}^{n}{n\choose k}^2{2k\choose k}\frac{1}{2k-1},
\end{align*}
and $S_n^{+}=\sum_{k=0}^{n}{n\choose k}^2{2k\choose k}(2k+1)^2$. In this paper we mainly prove that
\begin{align*}
4\sum_{k=0}^{n-1}kS_k\equiv \sum_{k=0}^{n-1}s_k\equiv \sum_{k=0}^{n-1}S_k^{+}\equiv 0\pmod{n^2}\quad\text{for $n\geqslant 1$},
\end{align*}
by establishing some binomial coefficient identities, such as
\begin{align*}
4\sum_{k=0}^{n-1}kS_k=n^2\sum_{k=0}^{n-1}\frac{1}{k+1}{2k\choose k}\left(6k{n-1\choose k}^2+{n-1\choose k}{n-1\choose k+1}\right).
\end{align*}
This confirms several recent conjectures of Z.-W. Sun.

\vskip 3mm \noindent {\it Keywords}: congruences; Legendre symbol; Zeilberger's algorithm;

\vskip 2mm
\noindent{\it MR Subject Classifications}: 11A07, 11B65, 05A10

\section{Introduction}
Inspired by the Schr\"{o}der numbers in combinatorics,
Z.-W. Sun \cite{Sun} introduced the following two kinds of numbers
\begin{align*}
R_n=\sum_{k=0}^{n}{n+k\choose 2k}{2k\choose k}\frac{1}{2k-1},\qquad
S_n=\sum_{k=0}^{n}{n\choose k}^2{2k\choose k}(2k+1),
\end{align*}
and proved many remarkable arithmetic properties of these numbers. For example, he proved that, for any odd prime $p$ and
positive integer $n$,
\begin{align}
\sum_{k=0}^{p-1}R_k &\equiv -p-\left(\frac{-1}{p}\right) \pmod{p^2}, \notag\\
\frac{1}{n^2}\sum_{k=0}^{n-1}S_k &=\sum_{k=0}^{n-1}{n-1\choose k}^2{2k\choose k}\frac{1}{k+1}, \label{eq:sun0}
\end{align}
where $\big(\frac{\cdot}{p}\big)$ denotes the Legendre symbol.

The aim of this paper is to prove the following results.

\begin{thm}\label{thm:sun1}
Let $n$ be a positive integer and $p$ a prime. Then
\begin{align}
4\sum_{k=0}^{n-1}kS_k &\equiv  0\pmod{n^2}, \label{eq:sun11}  \\[5pt]
\sum_{k=0}^{p-1}kS_k  &\equiv \frac{p^2}{8}\left(5-9\left(\frac{p}{3}\right)\right) \pmod{p^3}. \label{eq:sun12}
\end{align}
\end{thm}

Let
\begin{align*}
s_n=\sum_{k=0}^{n}{n\choose k}^2{2k\choose k}\frac{1}{2k-1},
\quad\text{and}\quad S_n^{+}&=\sum_{k=0}^{n}{n\choose k}^2{2k\choose k}(2k+1)^2.
\end{align*}

\begin{thm}\label{thm:sun2}
Let $n$ be a positive integer and $p$ a prime. Then
\begin{align}
\sum_{k=0}^{n-1} s_k &\equiv  0\pmod{n^2}, \label{eq:sun21} \\[5pt]
\sum_{k=0}^{p-1}s_k &\equiv -\frac{p^2}{2}\left(1+9\left(\frac{p}{3}\right)\right) \pmod{p^3}, \label{eq:sun22}
\end{align}
\end{thm}

\begin{thm}\label{thm:sun3}
Let $n$ be a positive integer and $p$ a prime. Then
\begin{align}
\sum_{k=0}^{n-1} S_k ^{+}&\equiv  0\pmod{n^2}, \label{eq:sun31}\\[5pt]
\sum_{k=0}^{p-1} S_k^{+} &\equiv -p^2\left(\frac{p}{3}\right) \pmod{p^3}. \label{eq:sun32}
\end{align}
\end{thm}

The congruences \eqref{eq:sun11}--\eqref{eq:sun21}, and \eqref{eq:sun31} were
originally conjectured by Z.-W. Sun \cite[Conjectures 5.5 and 5.6]{Sun}.
Note that, by establishing a lemma on sums of $q$-binomial coefficients, Z.-W. Sun himself could prove
the congruences \eqref{eq:sun21} and \eqref{eq:sun31} modulo $n$.

The key idea in this paper is to find out new expressions of the
sums $\frac{4}{n^2}\sum_{k=0}^{n-1}kS_k$,
$\frac{1}{n^2}\sum_{k=0}^{n-1}s_k$ and
$\frac{1}{n^2}\sum_{k=0}^{n-1}S_k^{+}$, just like Sun's identity
\eqref{eq:sun0}. No doubt that Zeilberger's algorithm (see
\cite{Koepf,PWZ}) will play an important role in our proof of
Theorems \ref{thm:sun1}--\ref{thm:sun3}. We shall give two proofs of
\eqref{eq:sun21}, one of which also gives the following new
congruence:
\begin{align}
\sum_{k=0}^{n-1}{n-1\choose k}^2{2k+1\choose k}\frac{3}{4k^2-1}\equiv 0\pmod{4n-1}.\label{eq:sun33}
\end{align}

\section{Proof of Theorem \ref{thm:sun1}}
\noindent{\it Proof of \eqref{eq:sun11}.} It is well known that $\frac{1}{k+1}{2k\choose k}$ is an integer (the $k$th Catalan number).
We shall prove the congruence \eqref{eq:sun11} by establishing the following identity:
\begin{align}
4\sum_{k=0}^{n-1}kS_k=n^2 f_{n-1}, \label{eq:s4}
\end{align}
where
\begin{align*}
f_n=\sum_{k=0}^{n}\frac{1}{k+1}{2k\choose k}\left(6k{n\choose k}^2+{n\choose k}{n\choose k+1}\right).
\end{align*}

It is clear that \eqref{eq:s4} holds for $n=1$. It remains to show that
\begin{align}
4nS_n=(n+1)^2f_n-n^2f_{n-1},  \label{eq:sn-fn}
\end{align}
for all positive integers $n$. Define $u_n=4nS_n$.
Applying Zeilberger's algorithm, we find that the numbers $u_n$ satisfy the following recurrence:
\begin{align}
&\hskip -2mm
n(n+1)(n+2)(n+3)u_{n+3}-n(n+1)(n+3)(11n+29)u_{n+2} \nonumber\\
&+n(n+2)(19{n}^{2}+74n+87)u_{n+1}-9(n+1)^{3}(n+2)u_n=0. \label{eq:rec1}
\end{align}
It is interesting that we can deduce the following second-order recurrence for $u_n$ from \eqref{eq:rec1}:
\begin{align}
&\hskip -2mm
n(n+1)(n+2)(4n+3)(4n+7)u_{n+2}-n(4n+3)(4n+11)(10{n}^{2}+30n+23)u_{n+1}  \nonumber \\
&+9(n+1)^{3}(4n+11)(4n+7)u_{n}=0.  \label{eq:s1}
\end{align}
In fact, if we denote the left-hand sides of \eqref{eq:rec1} and \eqref{eq:s1} by $\alpha_n$ and $\beta_n$, respectively,
then we can easily check that
\begin{align}
(4n+11)(4n+7)\alpha_n+n\beta_{n+1}-(n+2)\beta_n=0. \label{eq:s2}
\end{align}
Therefore, by induction on $n$, we immediately get $\beta_n=0$, i.e., the recurrence \eqref{eq:s1} is true.

Let $v_n=(n+1)^2f_n-n^2f_{n-1}$. Then Zeilberger's algorithm gives the following relation:
\begin{align}
&\hskip -2mm
(n+2)(n+3)(128{n}^{4}+864{n}^{3}+2016{n}^{2}+1994n+693)v_{n+3} \nonumber \\
&-(1408{n}^{6}+17696{n}^{5}+88512{n}^{4}+225582{n}^{3}+309049{n}^{2}+215886n+59535)v_{n+2}  \nonumber \\
&+(2432{n}^{6}+30880{n}^{5}+155712{n}^{4}+399646{n}^{3}+550013{n}^{2}+384657n+106920)v_{n+1}\nonumber \\
&-9(n+1)^{2}(128{n}^{4}+1376{n}^{3}+5376{n}^{2}+9130n+5695)v_n=0. \label{eq:rec2}
\end{align}
Similarly as before, we can deduce the following simpler recurrence for $v_n$ from \eqref{eq:rec2}:
\begin{align}
&\hskip -2mm
n(n+1)(n+2)(4n+3)(4n+7)v_{n+2}-n(4n+3)(4n+11)(10{n}^{2}+30n+23)v_{n+1}  \nonumber \\
&+9(n+1)^{3}(4n+11)(4n+7)v_{n}=0,  \label{eq:s1-new}
\end{align}
by noticing that
\begin{align*}
&(4n+11)(4n+7)(n+1)\gamma_n+(128{n}^{4}+864{n}^{3}+2016{n}^{2}+1994n+693)\delta_{n+1} \nonumber\\
&-(128{n}^{4}+1376{n}^{3}+5376{n}^{2}+9130n+5695)\delta_n=0, 
\end{align*}
where $\gamma_n$ and $\delta_n$ denote the left-hand sides of \eqref{eq:rec2} and \eqref{eq:s1-new}, respectively.

Thus, we have proved that the sequences $u_n$ and $v_n$ satisfy the same recurrence \eqref{eq:s1} (i.e., \eqref{eq:s1-new}).
Also, it is easy to verify that $u_n=v_n$ for $n=1,2$. This proves that $u_n=v_n$ for all positive integers. Namely, the identity \eqref{eq:sn-fn}
is true. This completes the proof.  \qed

\medskip
\noindent{\it Proof of \eqref{eq:sun12}.} Letting $n=p$ be a prime in \eqref{eq:s4}, we obtain
\begin{align}
\sum_{k=0}^{p-1}kS_k=\frac{p^2}{4}\sum_{k=0}^{p-1}\frac{1}{k+1}{2k\choose k}\left(6k{p-1\choose k}^2+{p-1\choose k}{p-1\choose k+1}\right).
\label{eq:s4-p}
\end{align}
It is easy to see that ${p-1\choose k}\equiv (-1)^k \pmod{p}$ for $0\leqslant k\leqslant p-1$. Therefore, from \eqref{eq:s4-p} we deduce that
\begin{align*}
\sum_{k=0}^{p-1}kS_k\equiv \frac{p^2}{4}
\left(\sum_{k=0}^{p-1}\frac{6k}{k+1}{2k\choose k}-\sum_{k=0}^{p-2}\frac{1}{k+1}{2k\choose k}\right)
\pmod{p^3}.
\end{align*}
The proof then follows from the following two congruences due to Pan and Sun \cite[(1.4) and (1.16) with $a=0$]{PS} (see also \cite{ST2}):
\begin{align}
\sum_{k=0}^{p-1}{2k\choose k}&\equiv \left(\frac{p}{3}\right)\pmod{p}, \label{eq:ps-1} \\[5pt]
\sum_{k=0}^{p-1}\frac{1}{k+1}{2k\choose k}&\equiv \frac{1}{2}\left(3\left(\frac{p}{3}\right)-1\right)\pmod{p},   \label{eq:ps-2}
\end{align}
and the fact that
\begin{align}
\frac{1}{p}{2p-2\choose p-1}\equiv -1\pmod{p}.   \label{eq:ps-3}
\end{align}

\section{Proof of Theorem \ref{thm:sun2}}
We shall give two proofs of \eqref{eq:sun21} by establishing two different identities.

\medskip
\noindent{\it First proof of \eqref{eq:sun21}.}
We want to show that
\begin{align}
\sum_{k=0}^{n-1}s_k=n^2 g_{n-1}, \label{eq:sun34}
\end{align}
where
$$
g_n=\sum_{k=0}^{n}\frac{1}{k+1}{2k\choose k}\left(2{n\choose k}{n\choose k+1}-{n\choose k}^2\right).
$$

Let $w_n=(n+1)^2g_n-n^2g_{n-1}.$
Applying Zeilberger's algorithm, we find that $s_n$ and $w_n$ satisfy the same recurrence:
\begin{align}
(n+3)^{2}s_{n+3}-(11{n}^{2}+46n+47)s_{n+2}
+(19{n}^{2}+58n+63)s_{n+1}-9(n+1)^{2}s_{n}=0. \label{eq:rec3}
\end{align}
Moreover, it is easy check that $s_n=w_n$ for $n=1,2,3$.
This proves that $s_n=w_n$ for all positive integers $n$. The identity \eqref{eq:sun34} then follows from
the fact that $s_0=g_0=-1$.
\qed

\medskip
\noindent{\it Second proof of \eqref{eq:sun21}.} We shall prove that
\begin{align}
\sum_{k=0}^{n-1}s_k=\frac{n^2}{4n-1}\sum_{k=0}^{n-1}{n-1\choose k}^2{2k+1\choose k}\frac{3}{4k^2-1}. \label{eq:sun35}
\end{align}

We claim that the numbers $s_n$ satisfy a second-order recurrence:
\begin{align}
(n+2)^{2}(8n+5)s_{n+2}-(80{n}^{3}+234{n}^{2}+214n+63)s_{n+1}+9(n+1)^{2}(8n+13)s_{n}=0.  \label{eq:rec4}
\end{align}
Denote the left-hand sides of \eqref{eq:rec3} and \eqref{eq:rec4} by $\lambda_n$ and $\mu_n$, respectively. Then it is easy
to verify that
\begin{align*}
(8n+13)\lambda_n+\mu_{n+1}-\mu_{n}=0,
\end{align*}
which leads to \eqref{eq:rec4} by induction on $n$.

Let $$h_n=\frac{1}{4n+3}\sum_{k=0}^{n}{n\choose k}^2{2k+1\choose k}\frac{3}{4k^2-1},$$
and $x_n=(n+1)^2h_n-n^2h_{n-1}$. Applying Zeilberger's algorithm for $x_n$, we obtain
\begin{align}
&\hskip -2mm
(n+3)^{2}(4n+15) (512{n}^{3}+2752{n}^{2}+4504n+2069)x_{n+3}-(22528{n}^{6}+297728{n}^{5} \nonumber\\
&+1584608{n}^{4}+4336564{n}^{3}+6403785{n}^{2}+4793214n+1395765)x_{n+2} \nonumber\\
&+(38912{n}^{6}+453376{n}^{5}+2138848{n}^{4}+5339252{n}^{3}+7595337{n}^{2}+5794002n  \nonumber\\
&+1757133)x_{n+1}-9(n+1)^{2}(4n-1)(512{n}^{3}+4288{n}^{2}+11544n+9837) x_n=0.  \label{eq:rec5}
\end{align}
Similarly as before, the numbers $x_n$ also satisfy a simpler recurrence:
\begin{align}
(n+2)^{2}(8n+5)x_{n+2}-(80{n}^{3}+234{n}^{2}+214n+63)x_{n+1}+9(n+1)^{2}(8n+13)x_{n}=0, \label{eq:rec6}
\end{align}
since
\begin{align*}
&\hskip -2mm
(8n+13)\nu_n-(4n+15) (512{n}^{3}+2752{n}^{2}+4504n+2069)\tau_{n+1}\\
&+(4n-1) (512{n}^{3}+4288{n}^{2}+11544n+9837)\tau_{n}=0,
\end{align*}
where $\nu_n$ and $\tau_n$ denote the left-hand sides of \eqref{eq:rec5} and \eqref{eq:rec6}, respectively.

Thus, we have proved that $s_n$ and $x_n$ satisfy the same recurrence \eqref{eq:rec4} (i.e., \eqref{eq:rec6}).
The proof of \eqref{eq:sun35} then follows from the fact that $s_1=x_1$ and $s_2=x_2$.

Note that $\gcd(n^2,4n-1)=1$ and
$$
{2k+1\choose k}\frac{3}{4k^2-1}={2k\choose k}\left(\frac{2}{2k-1}-\frac{1}{k+1}\right)
$$ is an integer. The identity \eqref{eq:sun35} immediately implies that both \eqref{eq:sun21} and \eqref{eq:sun33} hold.
\qed

\medskip
\noindent{\it Proof of \eqref{eq:sun22}.} Letting $n=p$ be a prime in \eqref{eq:sun34}, we get
\begin{align*}
\sum_{k=0}^{p-1}s_k
&=p^2\sum_{k=0}^{p-1}\frac{1}{k+1}{2k\choose k}\left(2{p-1\choose k}{p-1\choose k+1}-{p-1\choose k}^2\right) \\
&\equiv p^2\left(-\sum_{k=0}^{p-2}\frac{2}{k+1}{2k\choose k}-\sum_{k=0}^{p-1}\frac{1}{k+1}{2k\choose k}\right) \pmod{p^3}.
\end{align*}
The proof then follows from the congruences \eqref{eq:ps-2} and \eqref{eq:ps-3}.
\qed

\section{Proof of Theorem \ref{thm:sun3}}
\noindent{\it Proof of \eqref{eq:sun31}.}
We need to prove the following identity:
\begin{align}
\sum_{k=0}^{n-1}S_k^{+}=n^2\sum_{k=0}^{n-1}\left({2k+1\choose k}{n-1\choose k}^2
+\frac{1}{k+1}{2k\choose k}{n-1\choose k}{n-1\choose k+1}\right). \label{eq:s5}
\end{align}
Applying Zeilberger's algorithm, we find that the numbers $S_n^{+}$ satisfy the following recurrence:
\begin{align}
&\hskip -2mm
(n+3)^{2}(16{n}^{2}+69n+71)S_{n+3}^{+}-(176{n}^{4}+1911{n}^{3}+7600{n}^{2}+13101n+8217)S_{n+2}^{+} \notag\\
&+(304{n}^{4}+2943{n}^{3}+10721{n}^{2}+18171n+12672)S_{n+1}^{+} \notag\\
&-9(n+1)^{2}(16{n}^{2}+101n+156) S_{n}^{+}=0.  \label{eq:rec7}
\end{align}
It follows by induction on $n$ that
\begin{align}
&\hskip -2mm
(n+2)^{2}(256{n}^{4}+1280{n}^{3}+2176{n}^{2}+1488n+363)S_{n+2}^{+} \notag\\
&-(2560{n}^{6}+25600{n}^{5}+100096{n}^{4}+194848{n}^{3}+198238{n}^{2}+99610n+19677)S_{n+1}^{+} \notag\\
&+9(n+1)^{2}(256{n}^{4}+2304{n}^{3}+7552{n}^{2}+10704n+5563)S_{n}^{+}=0,  \label{eq:rec8}
\end{align}
since we have
\begin{align*}
&(256{n}^{4}+2304{n}^{3}+7552{n}^{2}+10704n+5563)\xi_n+(16{n}^{2}+69n+71)\eta_{n+1}\\
&-(16{n}^{2}+101n+156)\eta_{n}=0,
\end{align*}
where $\xi_n$ and $\eta_n$ denote the left-hand sides of \eqref{eq:rec7} and \eqref{eq:rec8}, respectively.

Denote the right-hand side of \eqref{eq:s5} by $n^2e_{n-1}$, and let $y_n=(n+1)^2e_{n}-n^2e_{n-1}$.
Applying Zeilberger's algorithm for $y_n$, we obtain
\begin{align}
&\hskip -2mm
 (n+3)^{2}(1024{n}^{5}+9984{n}^{4}+36736{n}^{3}+64000{n}^{2}+53236n+17057)y_{n+3}  \notag\\
&-(11264{n}^{7}+186624{n}^{6}+1275008{n}^{5}+4648064{n}^{4}+9750908{n}^{3}+11764759{n}^{2} \notag\\
&+7570338n+2011779)y_{n+2}+(19456{n}^{7}+324864{n}^{6}+2232448{n}^{5}+8173312{n}^{4} \notag\\
&+17192092{n}^{3}+20752931{n}^{2}+13332438n+3555639)y_{n+1} \notag\\
&-9(n+1)^{2}(1024{n}^{5}+15104{n}^{4}+86912{n}^{3}+244352{n}^{2}+336500n+182037)y_{n}=0. \label{eq:rec9}
\end{align}
It follows that the recurrence \eqref{eq:rec8} also holds for $y_n$, in view of
\begin{align*}
&(256{n}^{4}+2304{n}^{3}+7552{n}^{2}+10704n+5563)\zeta_n\\
&+(1024{n}^{5}+9984{n}^{4}+36736{n}^{3}+64000{n}^{2}+53236n+17057)\sigma_{n+1}\\
&-(1024{n}^{5}+15104{n}^{4}+86912{n}^{3}+244352{n}^{2}+336500n+182037)\sigma_{n}=0,
\end{align*}
where $\zeta_n$ and $\sigma_n$ denote the left-hand sides of \eqref{eq:rec9} and \eqref{eq:rec8} (with $S_n^{+}$ replaced by $y_n$ and so on),
respectively.  \eqref{eq:s5}. Noticing that $S_n^{+}=y_n$ for $n=1,2$ and $S_0^{+}=1$, we complete the proof.
\qed

\medskip
\noindent{\it Proof of \eqref{eq:sun32}.} Letting $n=p$ be a prime in \eqref{eq:s5}, we have
\begin{align*}
\sum_{k=0}^{p-1}S_k^{+}
&=p^2\sum_{k=0}^{p-1}\left({2k+1\choose k} {p-1\choose k}^2+\frac{1}{k+1}{2k\choose k}{p-1\choose k}{p-1\choose k+1}\right)\\
&\equiv p^2\left(\sum_{k=0}^{p-1}\frac{2k+1}{k+1}{2k\choose k}-\sum_{k=0}^{p-2}\frac{1}{k+1}{2k\choose k}\right) \pmod{p^3}.
\end{align*}
The proof then follows from the congruences \eqref{eq:ps-1}--\eqref{eq:ps-3}.
\qed

\section{Concluding remarks and open problems}
It is worth mentioning that Z.-W. Sun \cite[Conjecture 5.7]{Sun} also gave
an interesting $q$-analogue of \eqref{eq:sun21}. We hope that our proof of \eqref{eq:sun21} will give some hints to
tackle this $q$-analogue.

It seems that the congruences \eqref{eq:sun11} and \eqref{eq:sun31} can be further generalized as follows. Let
\begin{align*}
S_n^{(r)}=\sum_{k=0}^{n}{n\choose k}^2{2k\choose k}(2k+1)^r, \quad\text{and}\quad
T_n^{(r)}=\sum_{k=0}^{n}{n\choose k}^2{2k\choose k}(2k+1)^r (-1)^k.
\end{align*}
Numerical calculation suggests the following conjectures.
\begin{conj}
Let $n$ and $r$ be positive integers and $p$ a prime. Then
\begin{align}
\sum_{k=0}^{n-1}S_k^{(2r)} &\equiv  0\pmod{n^2}, \notag\\ 
\sum_{k=0}^{n-1}T_k^{(2r)} &\equiv  0\pmod{n^2}, \label{eq:conj-0}  \\
\sum_{k=0}^{p-1}T_k^{(2)} &\equiv  \frac{p^2}{2}\left(5-3\left(\frac{p}{5}\right)\right)\pmod{p^3}.
\end{align}
\end{conj}

Note that the $r=1$ case of \eqref{eq:conj-0} was conjectured by Z.-W. Sun (see \cite[Conjecture 5.6]{Sun}).

\begin{conj}Let $n$ and $r$ be positive integers. Then there exist integers $a_{2r-1}$ and $b_r$, independent of $n$, such that
\begin{align}
a_{2r-1}\sum_{k=0}^{n-1}S_k^{(2r-1)} &\equiv  0\pmod{n^2}, \label{eq:conj-2} \\
b_r\sum_{k=0}^{n-1}kS_k^{(r)} &\equiv  0\pmod{n^2}. \label{eq:conj-3}
\end{align}
\end{conj}

It seems rather difficult (almost impossible) to find out the best possible values of $a_{2r-1}$ and $b_r$ for all $r$.
Nevertheless, for small $r$, we propose the following conjecture.
\begin{conj}
Let $a_3=3$, $a_5=15$, $a_7=21$, $a_9=15$, $a_{11}=33$, $a_{13}=1365$, and $a_{15}=3$.
Let $b_2=12$, $b_3=4$, $b_4=60$, $b_5=20$, $b_6=84$, $b_7=28$, $b_8=60$, $b_9=20$, $b_{10}=132$,
$b_{11}=44$, and $b_{12}=5460$. Then the congruences \eqref{eq:conj-2} and \eqref{eq:conj-3} hold.
\end{conj}

\vskip 5mm \noindent{\bf Acknowledgments.} This work was partially
supported by the Fundamental Research Funds for the Central
Universities and the National Natural Science Foundation of China
(grant 11371144).

\end{document}